
\input amsppt.sty
\magnification=\magstep1
\TagsOnRight
\NoBlackBoxes
\nologo
\topmatter
\title GLEASON'S PROBLEM IN WEIGHTED BERGMAN SPACE ON EGG DOMAINS
\endtitle
\author Guangbin Ren and Jihuai Shi\endauthor
\affil   Department of Mathematics,\\
             University of Science and Technology of China,\\
            Hefei, Anhui, 230026,
            People's Republic of China\endaffil
\abstract
  In the paper, we discuss on the egg domains:
$$ \Omega_a=\left\{\xi=(z,w)\in\bold C^{n+m}: \ z\in\bold C^n, \
w\in\bold C^m, |z|^2+|w|^{2/a}<1\right\},      \qquad 0<a\le 2. $$
 We show that Gleason's problem can be solved in the weight Bergman
 space on the egg domains. The proof will need the help of the recent work of
 the second named  author on the weighted Bergman projections on this kind of
 domain. As an application, we obtain a multiplier theorem
 on the egg domains.
\endabstract

\keywords Gleason problem, Bergman space  \endkeywords
\footnote " " {1991 Mathematics Subject Classification : 32A10, 32A30}
\thanks Supported by the National Natural Sciences foundation of China
 and the National  Education Committee Doctoral foundation of China \endthanks

\endtopmatter
\document

{\heading  \S 1.  Introduction \endheading }
\par\medskip
  In the paper, we consider the egg domains:
 $$ \Omega_a=\left\{\xi=(z,w)\in\bold C^{n+m}: \ z\in\bold C^n, \
     w\in\bold C^m, |z|^2+|w|^{2/a}<1  \right\}, \qquad 0<a\le 2, $$
where $z=(z_1, \cdots, z_n)$, $w=(w_1, \cdots, w_m)$,
$|z|^2=\sum_{j=1}^{n} |z_j|^2$, $|w|^2=\sum_{j=1}^{m} |w_j|^2$. We also write
$\xi=(\xi_1, \cdots, \xi_{n+m})$.
If $0<a\le 2$,  then $\Omega_a$
is a pseudoconvex domain with  $C^1$ boundary and
$m$ pseudoconvex directions.
                                                        \par
For  $\xi\in\Omega_a$, we put
     $$     h(\xi)=h(z,w)=(1-|z|^2)^a-|w|^2,     $$
and for $\sigma>-1, \quad 1\le p<\infty$, let  $L^p(\Omega_a, dv_{\sigma})$
denote the space of measurable functions on $\Omega_a$ for which
  $$   \int_{\Omega_a} h^{\sigma}(\xi)|f(\xi)|^p dv(\xi)<\infty,  $$
where $dv$ is the volume measure on  $\Omega_a$,
$dv_{\lambda}=h^{\lambda}dv$.
                                                                \par
 As usual, $H(\Omega_a)$ is the space of all holomorphic functions on $\Omega_a $,
 $$A^p_{\sigma}(\Omega_a)=L^p(\Omega_a, dv_{\sigma})\cap H(\Omega_a)$$
denotes the weighted  Bergman space.
                                                                \par
Let $X$ be some class of holomorphic functions in a domain
 $\Omega\subset\bold C^N$. Gleason's problem, denoted as $(X, a, \Omega)$,
 is the following:
                                                                  \par
 For any given  $a\in\Omega$  , $f\in X$
and $f(a)=0$, Do there exist functions $f_1, \cdots, f_N\in X$, such that
 $f(z)=\sum_{k=1}^N (z_k-a_k)f_k(z)$?
                                                                \par
 The difficulty of the Gleason's problem depends on $\Omega$ and function space
 $X$. Gleason originally asked the problem for  $(B_n, 0, \text{ball algebra})$,
where $B_n$ is the unit ball of $C^n$.
This problem was solved by Leibenson. Subsequently, in the unit ball,
Rudin[5], Zhu[7], Ren and Shi[4] respectively discussed the following Gleason
problem: $(B_n, 0, C^k(B_n))$,
\quad $(B_n, 0, L^p_a(B_n))$, \quad $(B_n, 0, H_{p,q}(\varphi))$;
and in strongly pseudoconvex domain $\Omega$, Kerzman-Nagel [3], Ahern-Schneider [1]
 studied the Gleason's problem in Lipschitz space and $C^k$ space.
                                                                \par
  In this paper,
 We will prove that Gleason's problem
 $(\Omega_a, 0, A^p_{\lambda}(\Omega_a))
 \quad (1<p<\infty, \lambda\ge 0 \quad or \quad p=1,  \lambda>-1)$
 can be solved.  Its proof based on the recent work of
 the second named  author [6] on the weighted Bergman projections on the
 egg domains.  Our main result is the following:
                                                             \medskip

     \proclaim {Theorem A} Gleason's problem can be solved on $
A_{\lambda}^p(\Omega_a)
 \quad (1<p<\infty, \lambda\ge 0 \quad or \quad p=1,  \lambda>-1)$.
Furthermore,  for any
 $k\ge 1$, there exist bounded linear operators $A_{\alpha}$
$(|\alpha|=k)$
on $A_{\lambda}^p(\Omega_a)$, such that if
$ f\in A_{\lambda}^p(\Omega_a),  D^{\alpha}f(0)=0 \ (|\alpha|\le m-1) $, then $$
f(z)=\sum_{|\alpha|=m}z^{\alpha}A_{\alpha}f(z) $$ on $\Omega_a $, where
$ \alpha=(\alpha_1, \cdots,\alpha_n) $  be multiindex, $
|\alpha|=\alpha_1+\cdots+\alpha_n .$   \par  \endproclaim
                                                 \bigskip
As a direct corollary, we obtain a multiplier theorem.
\proclaim{Theorem B}
For $k=1,\cdots, n+m $, the transformation
$$
\sum_{\alpha} c_{\alpha}\xi^{\alpha}
\longrightarrow
\sum_{|\alpha|\ne 0}\frac{\alpha_k}{|\alpha|} c_{\alpha}\xi^{\alpha}
$$
maps $A^p_{\lambda}(\Omega_a)$
into
 $A^p_{\lambda}(\Omega_a)$.
 In other words, the complex sequence
 $\{\frac{\alpha_k}{|\alpha|} \}$
 is a multiplier of
  $A^p_{\lambda}(\Omega_a)$
  into
   $A^p_{\lambda}(\Omega_a)$.
\endproclaim

 {\heading \S 2.  Some Lemmas\endheading }
 \medskip
 \par
 For  $\sigma>-1$, let $K_{\sigma}$ be the Bergman Kernel function on
$A_{\sigma}^2(\Omega_a)$, then from [6],
 $$
 K_{\sigma}(\xi,\xi')
 =\sum_{k=0}^{n+1} c_k
 \frac{(1-<z,z'>)^{ak-n-1}}{((1-<z,z'>)^a-<w,w'>)^{\sigma+m+k}},       \tag 1
 $$
where $\xi=(z,w), \xi'=(z',w')$ are  points in $\Omega_a$, and, $c_k$
are constants only depending on
$m,n,\sigma,a$.
                                                \par
If $0<a\le 1$, denote
  $$
  G_{\sigma}(\xi,\xi')
  =\frac{(1-<z,z'>)^{(a-1)(n+2)/2}}
{((1-<z,z'>)^a-<w,w'>)^{(\sigma+m+n+2)/2}}.
  $$
If $1<a\le 2$, denote
  $$
  G_{\sigma}(\xi,\xi')
  =\frac{(1-<z,z'>)^{(a-1)(n+1)/2}}
{((1-<z,z'>)^a-<w,w'>)^{(\sigma+m+n+2)/2}}.
  $$
We will only discuss the case $0<a\le 1$, since the case
$1<a\le 2$ is similar.
 \par
 As usual, the symbol $A\lesssim B$ means that there exists a constant C
 such that $A\le CB$.

 \medskip
\par
\proclaim{Lemma 1 }
$$\left|\frac{\partial K_{\sigma}}{\partial \xi_k}
(\xi,\xi')\right|\lesssim |G_{\sigma}(\xi,\xi')|^2,   \qquad (1\le k \le n+m).
           \tag 2 $$
\endproclaim
                                                        \par
\demo{Proof}
 If $\xi=(z,w)$,  $\xi'=(z',w')\in\Omega_a$, then
    $$
    |z|^2+|w|^{2/a}<1,  \quad
    |z'|^2+|w'|^{2/a}<1.
     $$
Namely
$$\aligned
  |1-<z,z'>|^{2a}&\ge (1-|<z,z'>|)^{2a}\ge (1-|z||z'|)^{2a}  \\
       &\ge(1-|z|^2)^a(1-|z'|^2)^a>|w|^2|w'|^2,
  \endaligned
$$
thus
$$
|<w,w'>|<|1-<z,z'>|^a.          \tag 3
$$
Differentiating the both sides of the formula in (1), we obtain by (3)
  $$
 \left|\frac{\partial K_{\sigma}}{\partial z_k}\right|\lesssim \sum_{j=0}^{n+1}
\left|\frac{(1-<z,z'>)^{aj-n-1+a-1}}{((1-<z,z'>)^a-<w,w'>)^{\sigma+m+j+1}}\right|.
  $$
Again by (3), the  $j$-th  $(\forall 0\le j\le n+1)$
 summand can be controlled by the $j+1$-th summand,
therefore, can be controlled by the
 $n+1$-th summand. This means (2) holds for any $1\le k\le n$.
                                                        \par
Similarly
$$\aligned
\left|\frac{\partial K_{\sigma}}{\partial w_k}\right|&\lesssim
\left|\frac{(1-<z,z'>)^{(a-1)(n+1)}}{((1-<z,z'>)^a-<w,w'>)^{\sigma+m+n+2}}\right|
                  \\
  &\le
\left|\frac{(1-<z,z'>)^{(a-1)(n+2)}}{((1-<z,z'>)^a-<w,w'>)^{\sigma+m+n+2}}\right|
                  \\
   &=|G_{\sigma}(\xi,\xi')|^2.
\endaligned
$$
Here we use the condition of $0<a\le 1$. This proves Lemma 1.
\enddemo
                               \medskip
                                 \par
     Put
       $$
       \psi_{k,r}(\xi,\xi')
       =\frac{<w,w'>^k}{(1-<z,z'>)^r}    \quad(k\in\bold N\cup \{0 \}, \quad  r>0).
       $$
Shi [6] proved that for $s>-1$,
  $$\aligned
 &\quad \int_{\Omega_a}
  h^s(\xi')|\psi_{k,r}(\xi,\xi')|^2dv(\xi')
  \\
 & =\frac{\pi^{n+m}k!\Gamma(s+1)\Gamma(a(s+k+m)+1)}
{\Gamma^2(r)\Gamma(s+k+m+1)}|w|^{2k}
  \sum_{j=0}^{\infty}
  \frac{\Gamma^2(j+r)|z|^{2j}}{\Gamma(a(s+k+m)+j+n+1)j!}.
  \endaligned
  $$
                            \medskip
                              \par
\proclaim{Lemma 2 }
             If $ 0<d<\sigma+1$, then
 $$   \int_0^1\int_{\Omega_a} h^{\sigma-d}(\xi')
   |G_{\sigma}(t\xi,\xi')|^2dv(\xi')dt
   \lesssim h^{-d}(\xi).       \tag 4
   $$
\endproclaim
\demo{Proof}
 By the formula
        $$
        \frac{1}{(1-t)^s}=\sum_{k=0}^{\infty}
\frac{\Gamma(k+s)}{k!\Gamma(s)}t^k   \qquad  (|t|<1, \quad s>0),
        $$
we obtain
$$\aligned
G_{\sigma}(\xi,\xi')&=
\frac{1}{(1-<z,z'>)^{\frac{a(\sigma+m)+n+2}{2}}}
\frac{1}{(1-\frac{<w,w'>}{(1-<z,z'>)^a})^{\frac{\sigma+m+n+2}{2}}}
                                        \\
   &=\sum_{k=0}^{\infty}\frac{\Gamma(k+c)}{k!\Gamma(c)}
     \frac{<w,w'>^k}{(1-<z,z'>)^{ak+b}}
                                                \\
   &= \sum_{k=0}^{\infty}\frac{\Gamma(k+c)}{k!\Gamma(c)}
      \psi_{k,ak+b}(\xi,\xi')
\endaligned
$$
where  $b=\frac{a(\sigma+m)+n+2}{2}$, $c=\frac{\sigma+m+n+2}{2}$.
                                                                \par
Denote  $\mu=\sigma-d$. Since $\{\psi_{k,ak+b}\}_{k=0}^{\infty}$ is a
 orthogonal basis on
$\Omega_a$ [6],  Parseval equality tells us  that
                                                                \par
  $$
  \aligned
  \text{Left side of (4)}
 &=\int_0^1\int_{\Omega_a}h^{\mu}(\xi')|G_{\sigma}(t\xi,\xi')|^2dv(\xi')dt
                                                                \\
 &=\sum_{l=0}^{\infty}\frac{\Gamma^2(l+c)}{(l!)^2\Gamma^2(c)}
 \int_0^1\int_{\Omega_a}h^{\mu}(\xi')|\psi_{l,al+b}(t\xi,\xi')|^2dv(\xi')dt
                                                                \\
&=\sum_{l=0}^{\infty}\frac{\Gamma^2(l+c)}{(l!)^2\Gamma^2(c)}
\frac{\pi^{n+m}l!\Gamma(\mu+1)\Gamma(a(\mu+l+m)+1)}
{\Gamma^2(al+b)\Gamma(\mu+l+m+1)}|w|^{2l}
                                                                \\
&\qquad \sum_{j=0}^{\infty}
\frac{\Gamma^2(j+la+b)|z|^{2j}}{\Gamma(a(\mu+l+m)+j+n+1)j!}
\frac{1}{2(l+j)+1};
 \endaligned
 $$
 $$
 \aligned
 \text{Right side of (4)}&=\frac{1}{((1-|z|^2)^a-|w|^2)^d}
                                        \\
   &=\frac{1}{(1-|z|^2)^{ad}}
\frac{1}{(1-\frac{|w|^2}{(1-|z|^2)^a})^d}
                                        \\
&=\sum_{l=0}^{\infty}\frac{\Gamma(l+d)}{\Gamma(d)l!}\sum_{j=0}^{\infty}
\frac{\Gamma(j+a(l+d))}{\Gamma(a(l+d))j!}|z|^{2j}|w|^{2l}.
 \endaligned
 $$
Comparing the coefficient of $|z|^{2j}|w|^{2l}$,
we only need to prove the following inequality
$$\aligned
&\frac{1}{2(l+j)+1}
\frac{\Gamma^2(l+c)\Gamma(a(\mu+l+m)+1)\Gamma^2(j+al+b)}
{\Gamma^2(la+b)\Gamma(\mu+l+m+1)\Gamma(a(\mu+l+m)+j+n+1)}
\\
&\lesssim\frac{\Gamma(l+d)\Gamma(j+a(l+d))}{\Gamma(a(l+d))},
\endaligned
$$
while this can be changed to prove
$$
\aligned
\frac{\Gamma^2(j+la+b)}
{(2(l+j)+1)\Gamma(a(\mu+l+m)+j+n+1)\Gamma(j+a(l+d))}
\lesssim 1,
\endaligned        \tag 5
$$
$$
\aligned
\frac{\Gamma^2(l+c)}
{\Gamma(l+d)\Gamma(\mu+l+m+1)}
\frac{\Gamma(a(\mu+l+m)+1)\Gamma(a(l+d))}{\Gamma^2(la+b)}
\lesssim 1,
\endaligned  \tag 6
$$
where  $b=\frac{a(\sigma+m)+n+2}{2}$, $c=\frac{\sigma+m+n+2}{2}$,
  $\mu=\sigma-d$.
                                                                \par
Thus it remains to prove
(5) and (6) for bigger enough
 $j, l$.
 This is not too hard to prove
 by the well known properties of  $\Gamma$ function,
 $$\aligned
   &(\text{i}). \qquad     \Gamma(x+1)=x\Gamma(x) \qquad (x>0),
   \\
   &(\text {ii}). \qquad     \frac{\Gamma^2(x+t)}{\Gamma(t)\Gamma(2x+t)}\le1
                     \qquad (x\ge 0, \quad t>0),
   \\
   &(\text {iii}). \qquad \lim_{t\longrightarrow\infty}
                  \frac{\Gamma^2(x+t)}{\Gamma(t)\Gamma(2x+t)}=1 \qquad (x\ge 0).
\endaligned
$$
In fact,
    $$
  \aligned
  \text{Left side of (5)}
 &\lesssim\frac{\Gamma^2(j+la+b)}
   {\Gamma(a(\mu+l+m)+j+n+1)\Gamma(j+a(l+d)+1)}
   \\
 &=
   \frac{\Gamma^2(j+la+\frac{a(\sigma+m)+n+2}{2})}
   {\Gamma(a(\sigma-d+l+m)+j+n+1)\Gamma(j+a(l+d)+1)}
   \\
   &\le 1,
 \endaligned
 $$
 where the last step follows from (ii) in the properties of $\Gamma$
 function, since in the fraction of the last inequality
 the sum of the two numbers
 in the brackets
 of  the denominator
 doubles the number in the bracket of the numerator;
 similarly we obtain
 $$
 \aligned
&\text{Left side of (6)}\\
&\lesssim
 \frac{\Gamma^2(l+c)}
 {\Gamma(l+d+n+1)\Gamma(\mu+l+m+1)}
 \frac{\Gamma(a(\mu+l+m)+1)\Gamma(a(l+d)+n+1)}{\Gamma^2(la+b)}
 \\
&=
\frac{\Gamma^2(l+\frac{\sigma+m+n+2}{2})}
 {\Gamma(l+d+n+1)\Gamma(\sigma-d+l+m+1)}
 \frac{\Gamma(a(\sigma-d+l+m)+1)\Gamma(a(l+d)+n+1)}
{\Gamma^2(la+\frac{a(\sigma+m)+n+2}{2})}
\\
&\lesssim 1.
 \endaligned
 $$
 This completes the proof of Lemma 2.
\enddemo
                                                         \bigskip
                                                         \par
{\heading  \S 3.  Proof of the theorems  \endheading }
\par
 To prove the main theorem, we also need the following known results.
\par
Let $K_{\sigma}$ be weighted Bergman Kernel function on the space
 $A_{\sigma}^2(\Omega_a)$, define the operator
$$
(T_{\sigma}f)(\xi)=
C_{\sigma}\int_{\Omega_a}h^{\sigma}(\xi')K_{\sigma}(\xi,\xi')f(\xi')dv(\xi'),
$$
where
 $
 C_{\sigma}=\left\{
 \int_{\Omega_a} h^{\sigma}(\xi')K_{\sigma}(\xi,\xi')dv(\xi')
 \right\}^{-1},
 $
 this is a constant independent of $\xi$ [6].
                                                        \medskip
                                                        \par
\proclaim{Lemma 3 ([6]) }
 For $1\le p<\infty$, $T_{\sigma}$ is a linear bounded  operator on
 $L^p(\Omega_a,dv_{\lambda})$
iff
$$\aligned
0<\lambda+1<p(\sigma+1)
\endaligned       \tag 7
$$
and when (7) holds, $T_{\sigma}$ is a bounded projection operator from
 $L^p(\Omega_a,dv_{\lambda})$
to $A^p_{\lambda}(\Omega_a)$.
\endproclaim
                                                \medskip
                                                \par
\proclaim{Lemma 4 (Schur Lemma [5]) }
       Let $(X,\mu)$ be a measurable space. Suppose $Q$ is a
       non-negative measurable function on
 $X\times X$, $1<p<\infty$, $\frac{1}{p}+\frac{1}{q}=1$.
For the  integral operator $T$ induced by $Q$, that is,
$Tf(x)=\int_X Q(x,y)f(y)d\mu(y)$,
if there exists a non-negative measurable function $g$ on $X$ and constant $C$
such that
$$
\int_X Q(x,y)g^q(y)d\mu(y)\le C g^q(x),  \qquad a.e. \qquad x\in X,       \tag 8
$$
$$
\int_X Q(x,y)g^p(x)d\mu(x)\le C g^p(y),  \qquad a.e. \qquad y\in X,       \tag 9
$$
then $T$ is the bounded operator on $L^p(X,d\mu)$, and $||T||\le C$.
\endproclaim
                                                     \bigskip
                                                     \par

Now we set to give the proof of the Main Theorem.
                                                 \par\medskip
\demo{\bf Proof of Theorem A}
 We only need to prove in the case $|\alpha|=1$, the general
 case can be proved by induction as in Zhu [7].
                                                                \par
 For $0<a\le 2$,  $\Omega_a$ is a convex Reinhardt domain.
 By the
 Leibenzon decomposition [5] on convex domain,
 for $f\in H(\Omega_a), f(0)=0,$
 $$
 f(\xi)=\sum_{k=1}^{n+m}\xi_k\int_0^1\frac{\partial f}{\partial \xi_k}(t\xi)dt.
 $$
Denote $$T_kf(\xi)=\int_0^1\frac{\partial f}{\partial\xi_k}(t\xi)dt.$$
Clearly we only need to demonstrated that $T_k$
is a bounded operator in
 $A^p_{\lambda}(\Omega_a)$ for each $k$.
                                                                        \par
Due to  Lemma 3, for $f\in A^p_{\lambda}(\Omega_a)$, there exists
the reproducing formula:
  $$
  f(\xi)=C_{\sigma}\int_{\Omega_a}h^{\sigma}(\xi')K_{\sigma}(\xi,\xi')f(\xi')dv(\xi')
  \tag 10
  $$
where
$$ 0<\lambda+1<p(\sigma+1).    \tag 11 $$
Hence
$$
\aligned
T_kf(\xi)&=\int_0^1\frac{\partial
f}{\partial \xi_k}(t\xi)dt
                                                \\
&=C_{\sigma}\int_0^1\int_{\Omega_a} h^{\sigma}(\xi')
\frac{\partial K_{\sigma}(t\xi,\xi')}{\partial\xi_k}f(\xi')dv(\xi')dt
                                                \\
&=C_{\sigma}\int_{\Omega_a}f(\xi')\left(
\int_0^1 h^{\sigma}(\xi')
\frac{\partial K_{\sigma}(t\xi,\xi')}{\partial\xi_k}dt
\right)dv(\xi')
                                                \\
&=\int_{\Omega_a}f(\xi')Q(\xi,\xi')dv(\xi').
\endaligned                     \tag12
$$
Here
$$ Q(\xi,\xi')=
C_{\sigma}\int_0^1 h^{\sigma}(\xi')
\frac{\partial K_{\sigma}(t\xi,\xi')}{\partial\xi_k}dt.
$$
So
 $T_k$ is just the integral operator induced by  $Q$.
                                                   \par
To prove the boundedness of $T_k$ in $A_{\lambda}^p(\Omega_a)$,
we treat the two cases $1\le p<\infty$ and $p=1$ separately.

                                                \par
 Case 1.  $1<p<\infty,  \lambda\ge 0$.
                                                  \par
    Choose $\sigma$, such that $\sigma>\frac{\lambda+1}{p}-1$,
    i.e. (11) holds.  Since
    $\lambda \ge 0$,  then $\sigma>\frac{1}{p}-1$, so the intersection
    of the two intervals
 $ (0, \sigma+1)\cap (-\frac{\sigma}{p-1}, \frac{1}{p-1})$ is nonempity.
 Pick
           $$
           d\in (0, \sigma+1)\cap (-\frac{\sigma}{p-1}, \frac{1}{p-1}).
           $$
    \par
   In the  Schur Lemma,
 take $X=\Omega_a, \quad d\mu=dv_{\lambda}, \quad g=h^{-\frac{d}{q}}$.
  Then  (8) and (9) turn into
$$
\int_{\Omega_a} |Q(\xi,\xi')|h^{-d+\lambda}(\xi')dv(\xi')\lesssim h^{-d}(\xi),
     \tag 13
$$
$$
\int_{\Omega_a} |Q(\xi,\xi')|h^{-d(p-1)+\lambda}(\xi)dv(\xi)
\lesssim h^{-d(p-1)}(\xi').
     \tag 14
$$
Note that $\lambda\ge 0$,  $h^{\lambda}\lesssim 1$, thus if
 (13) and (14) hold for $\lambda=0$,  they must be hold for any
 $\lambda\ge 0$.
On the other hand
$$
\aligned
|Q(\xi,\xi')|&=\left|C_{\sigma}\int_0^1 h^{\sigma}(\xi')
\frac{\partial K_{\sigma}(t\xi,\xi')}{\partial \xi_k}dt\right|
                                         \\
&\lesssim\int_0^1 h^{\sigma}(\xi')|G_{\sigma}(t\xi,\xi')|^2dt.
\endaligned              \tag 15
$$
Then it remains to prove:
$$
  \int_0^1\int_{\Omega_a} h^{\sigma-d}(\xi')
  |G_{\sigma}(t\xi,\xi')|^2dv(\xi')dt
  \lesssim h^{-d}(\xi),
  $$
$$
    \int_0^1\int_{\Omega_a} h^{-d(p-1)}(\xi)
    |G_{\sigma}(t\xi,\xi')|^2dv(\xi)dt
    \lesssim h^{-d(p-1)-\sigma}(\xi').       \tag 16
  $$
By Lemma 2, it only needs to show  that (16) holds.
                                                \par
By the definition of $G_{\sigma}$,
        $$
        |G_{\sigma}(t\xi,\xi')|=|G_{\sigma}(t\xi',\xi)|.
        $$
Again use Lemma 2,
     $$
     \aligned
    &\int_0^1\int_{\Omega_a} h^{-d(p-1)}(\xi)
    |G_{\sigma}(t\xi,\xi')|^2dv(\xi)dt
                                   \\
        & =   \int_0^1\int_{\Omega_a} h^{-d(p-1)}(\xi)
         |G_{\sigma}(t\xi',\xi)|^2dv(\xi)dt
                                        \\
        & \lesssim h^{-d(p-1)-\sigma}(\xi').
 \endaligned
     $$
     \par
     Case 2: $p=1, \lambda>-1$.
     \par
       Choose $\sigma$, such that $\sigma>\lambda>-1$, then (11) holds.
 As a consequence of Lemma  2, (12) and  (15) , we have
 $$ \aligned
 \int_{\Omega_a}|Q(\xi,\xi')|dv_{\lambda}(\xi)&\lesssim
 \int_{\Omega_a}\int_0^1 h^{\sigma}(\xi')|G_{\sigma}(t\xi,\xi')|^2dt
dv_{\lambda}(\xi)                                \\
&=h^{\sigma}(\xi')
\int_0^1 \int_{\Omega_a} |G_{\sigma}(t\xi,\xi')|^2h^{\lambda}(\xi)
dv(\xi)dt                                                \\
&\lesssim h^{\sigma}(\xi')h^{\lambda-\sigma}(\xi')=h^{\lambda}(\xi').
 \endaligned  $$

Therefore
$$\aligned
\int_{\Omega_a}|T_kf(\xi)|dv_{\lambda}(\xi)
&\le \int_{\Omega_a}\int_{\Omega_a}|f(\xi')|
|Q(\xi,\xi')|dv(\xi')dv_{\lambda}(\xi)            \\
&=\int_{\Omega_a}|f(\xi')|\left(\int_{\Omega_a}
|Q(\xi,\xi')|dv_{\lambda}(\xi)
\right)dv(\xi')                         \\
&\lesssim\int_{\Omega_a}|f(\xi')|dv_{\lambda}(\xi') .
\endaligned $$
This completes the proof of  Theorem A.
\enddemo
                                                          \par \medskip
\demo{\bf Proof of theorem B}
Since $\Omega_a$ is a Reinhardt domain, every holomorphic function on it
has Taylor expansion.
By the proof of Theorem A,
$f(\xi)=\sum_{\alpha}c_{\alpha}z^{\alpha}\in A_{\lambda}^p(\Omega_a)$ implies
$$\xi_k T_k f(\xi)=\xi_k \int_0^1\frac{\partial f}{\partial\xi_k}(t\xi)dt
=\sum_{|\alpha|\ne 0}\frac{\alpha_k}{|\alpha|}c_{\alpha}\xi^{\alpha}
\in A_{\lambda}^p(\Omega_a)$$
This completes the proof of Theorem B.

\enddemo

      \end{document}